\documentclass[aos,MSNbibl,seceqn,nameyear,dvips]{arximspdf}
\usepackage{algorithm}
\usepackage{algorithmic}
\usepackage{graphicx}

\doi{10.1214/12-AOS997} 
\volume{40}
\issue{2}
\pubyear{2012}
\firstpage{1074}
\lastpage{1101}

\makeatletter
\newproclaim{definition}{Definition}
\newtheorem{theorem}{Theorem}
\newproclaim{remark}{Remark}

\newcommand{\sign}{\operatorname{sgn}}

\renewcommand{\d}{\delta}
\newcommand{\s}{\sigma}
\newcommand{\T}{\Theta}

\renewcommand{\aa}{\mathcal{A}}
\newcommand{\bb}{\mathcal{B}}
\newcommand{\hh}{\mathcal{H}}

\renewcommand{\b}{\boldsymbol{\beta}}
\newcommand{\bhat}{\hat{\b}}
\newcommand{\corr}{\operatorname{corr}}
\newcommand{\eb}{\boldsymbol{\eps}}
\newcommand{\eboost}{elasticBoost}
\newcommand{\enet}{\mathrm{enet}}
\newcommand{\eps}{{\varepsilon}}
\newcommand{\floor}{\operatorname{floor}}
\newcommand{\F}{\mathbf{F}}
\newcommand{\Fhat}{\hat{\F}}
\newcommand{\FSe}{$\mathrm{FS}_\varepsilon$}
\newcommand{\g}{\mathbf{g}}
\newcommand{\I}{\mathbf{I}}
\newcommand{\lstar}{l^*}
\newcommand{\Ms}{{M^*}}
\newcommand{\nuhat}{\hat{\nu}}
\newcommand{\one}{\mathbf{1}}
\newcommand{\Rhat}{{\hat{R}}}
\newcommand{\rhohat}{\hat{\rho}}
\newcommand{\X}{\mathbf{X}}
\newcommand{\Xs}{\X^*}
\newcommand{\y}{\mathbf{y}}
\newcommand{\ys}{{\y^*}}
\newcommand{\Z}{\mathbf{Z}}
\newcommand{\zero}{\mathbf{0}}

\newcommand{\E}{\mathbb{E}}
\newcommand{\RR}{\mathbb{R}}
\def\mref#1{(\ref{#1})}
\def\l2boost{$L_2$Boosting}

\newcommand{\argmin}{{\arg\min}}%
\newcommand{\argmax}{{\arg\max}}%
\makeatother

\begin{document}
\begin{frontmatter}

\title{Characterizing $L_2$Boosting}
\runtitle{Characterizing $L_2$Boosting}

\begin{aug}
\author[A]{\fnms{John} \snm{Ehrlinger}\corref{}\ead[label=e1]{john.ehrlinger@gmail.com}}
\and
\author[B]{\fnms{Hemant} \snm{Ishwaran}\ead[label=e2]{hemant.ishwaran@gmail.com}\ead[label=u1,url]{http://web.ccs.miami.edu/\textasciitilde hishwaran}}
\runauthor{J. Ehrlinger and H. Ishwaran}
\affiliation{Cleveland Clinic and University of Miami}
\address[A]{Department of Quantitative Health Sciences\\
Cleveland Clinic\\
Cleveland, Ohio 44195\\
USA\\
\printead{e1}} 
\address[B]{Division of Biostatistics\\
Department of Epidemiology\\
\quad and Public Health\\
University of Miami\\
Miami, Florida 33136\\
USA\\
\printead{e2}\\
\printead{u1}}
\end{aug}

\received{\smonth{10} \syear{2011}}
\revised{\smonth{3} \syear{2012}}

%
\begin{abstract}
We consider $L_2$Boosting, a special case of Friedman's generic boosting
algorithm applied to linear regression under $L_2$-loss. We study
$L_2$Boosting for an arbitrary regularization parameter and derive an
exact closed form expression for the number of steps taken along a
fixed coordinate direction. This relationship is used to describe
$L_2$Boosting's solution path, to describe new tools for studying its
path, and to characterize some of the algorithm's unique properties,
including active set cycling, a property where the algorithm spends
lengthy periods of time cycling between the same coordinates when
the regularization parameter is arbitrarily small. Our fixed
descent analysis also reveals a \textit{repressible condition} that
limits the effectiveness of $L_2$Boosting in correlated problems by
preventing desirable variables from entering the solution path. As
a simple remedy, a data augmentation method similar to that used for
the elastic net is used to introduce $L_2$-penalization and is
shown, in combination with decorrelation, to reverse the repressible
condition and circumvents $L_2$Boosting's deficiencies in correlated
problems. In itself, this presents a new explanation for why the
elastic net is successful in correlated problems and why methods like
LAR and lasso can perform poorly in such settings.
\end{abstract}

%

\vspace*{6pt}
\end{frontmatter}

\section{Introduction}\label{S:introduction}

Given data $\{y_i, \mathbf{x}_i\}_{1}^n$, where $y_i$ is the response and
$\mathbf{x}_i=(x_{i,1},\ldots,x_{i,p})\in\RR^p$ is the $p$-dimensional
covariate, the goal in many analyses is to approximate the unknown
function $F(\mathbf{x})=\E(y|\mathbf{x})$ by minimizing a~specified loss function
$L(y, F)$ [a common choice is $L_2$-loss, $L(y,F)=(y-F)^2/2$]. In
trying to estimate $F$, one strategy is to make use of a large system
of possibly redundant functions $\hh$. If $\hh$ is rich enough, then
it is reasonable to expect $F$ to be well approximated by an additive
expansion of the form
\[
F(\mathbf{x}; \{\beta_k,\alpha_k\}_1^K) = \sum_{k=1}^K \beta_k h(\mathbf{x};\alpha_k),\vadjust{\goodbreak}
\]
where $h(\mathbf{x};\alpha)\in\hh$ are base learners parameterized by
$\alpha\in\Theta$. To estimate $F$, a~joint multivariable optimization
over $\{\beta_k,\alpha_k\}_1^K$ may be used. But such an optimization may
be computationally slow or even infeasible for large dictionaries.
Overfitting may also result. To circumvent this problem, iterative
descent algorithms are often used.

One popular method is the gradient descent algorithm described
by \citet{Friedman:2001}, closely related to the method of ``matching
pursuit'' used in the signal processing literature [\citet{MZ:1993}].
This algorithm is applicable to a wide range of problems and loss
functions, and is now widely perceived to be a generic form of
boosting. For the $m$th step, $m=1,\ldots,M$, one solves
%
\begin{equation}\label{friedman.line.search}
\rho_m = \mathop\argmin_{\rho\in\RR}\sum_{i=1}^n
L \bigl(y_i, F_{m-1}(\mathbf{x}_i) + \rho h(\mathbf{x}_i;\alpha_m) \bigr),
\end{equation}
where
%
\begin{equation}\label{friedman.grad.optimization}
\alpha_m=\mathop\argmin_{\alpha\in\Theta}\sum_{i=1}^n [g_m(\mathbf{x}_i)
- h(\mathbf{x}_i;\alpha) ]^2
\end{equation}
identifies the closest base learner to the gradient $\g_m =
(g_m(\mathbf{x}_1),\ldots, g_m(\mathbf{x}_n))^T$ in $L_2$-distance, where
$g_m(\mathbf{x}_i)$
is the gradient evaluated at the current value $F_{m-1}(\mathbf{x}_i)$,
and is
defined by
\[
g_m(\mathbf{x}_i)
=-\biggl [\frac{\partial L(y_i,F(\mathbf{x}_i))}{\partial F(\mathbf{x}_i)}
\biggr]_{F_{m-1}(\mathbf{x}_i)}
= -L'(y_i, F_{m-1}(\mathbf{x}_i)).
\]
The $m$th update for the predictor of $F$ is
\[
F_m(\mathbf{x}) = F_{m-1}(\mathbf{x}) + \nu\rho_m h(\mathbf{x};\alpha_m),
\]
where $0<\nu\le1$ is a regularization (learning) parameter.

In this paper, we study Friedman's algorithm under $L_2$-loss in
linear regression settings assuming an $n\times p$ design matrix
$\X=[\X_1,\ldots,\X_p]$, where $\X_k=(x_{1,k},\ldots,x_{n,k})^T$
denotes the $k$th column. Here $\X_k$ represents the $k$th
base learner; that is, $h(\mathbf{x}_i;k)=x_{i,k}$ where $k=\alpha$ and
$\T=\{1,\ldots,p\}$. It is well known that under $L_2$-loss the gradient
simplifies to the residual $g_m(\mathbf{x}_i) = y_i -
F_{m-1}(\mathbf{x}_i)$. This is particularly attractive for a theoretical
treatment as it allows one to combine the
line-search \mref{friedman.line.search} and the
learner-search~\mref{friedman.grad.optimization} into a single step
because the $L_2$-loss function can be expressed as $L(y_i,
F_{m-1}(\mathbf{x}_i) + \rho x_{i,k}) = (g_m(\mathbf{x}_i) - \rho
x_{i,k})^2$. The
optimization problem becomes
\[
\{\rho_m,k_m\} =
\mathop{\argmin}_{{\rho\in\RR, 1\le k\le p}}\| \g_m - \rho\X_k
\|^2.\vadjust{\goodbreak}
\]

It is common practice to standardize the response by removing its mean
which eliminates the issue of whether an intercept should be included
as a~column of~$\X$. It is also common to standardize the columns of
$\X$ to have a~mean of zero and squared-length of one.
Thus, throughout, we assume the data is
standardized according to
%
\begin{equation}
\label{data.standardization}\qquad
\sum_{i=1}^n y_i = 0,\qquad
\sum_{i=1}^n x_{i,k} = 0,\qquad
\sum_{i=1}^n x_{i,k}^2 = 1, \qquad k=1,\ldots,p.
\end{equation}
The condition $\sum_{i=1}^n x_{i,k}^2 = 1$ leads to a particularly useful
simplification:
\[
\rho_m = \X_{k_m}^T\g_m,\qquad
k_m = \mathop\argmax_{1\le k\le p}|\X_k^T\g_m|.
\]
Thus, the search for the most favorable direction is equivalent to
determining the largest absolute value $|\X_k^T\g_m|$. We refer to
$\X_k^T\g_m$ as the \textit{gradient-correlation} for $k$. We shall
refer to Friedman's algorithm under the above settings as \l2boost.
Algorithm~\ref{A:L2Boost} provides a formal description of the algorithm
[we use $\F_{m-1}=(F_{m-1}(\mathbf{x}_1),\ldots,F_{m-1}(\mathbf{x}_n))^T$
for notational convenience].

\begin{algorithm}[t]
\caption{\l2boost}\label{A:L2Boost}
\begin{algorithmic}[1]
\STATE Initialize $F_{0,i}=0$ for $i=1,\ldots,n$
\FOR{$m=1$ to $M$}
\STATE$k_m = \argmax_{1\le k\le p} |\X_k^T\g_m|$, where $\g_m=\y
-\F_{m-1}$
\STATE$\F_m = \F_{m-1} + \nu\rho_m\X_{k_m}$, where $\rho_m=\X
_{k_m}^T\g_m$
\ENDFOR
\end{algorithmic}
\end{algorithm}

Properties of stagewise algorithms similar to \l2boost have been
studied extensively under the assumption of an infinitesimally
small regularization parameter. \citet{EfronLAR} considered a forward
stagewise algorithm \FSe, and showed under a convex cone condition
that the Least Angle Regression (LAR) algorithm yields the solution
path for $\mathrm{FS}_{0}$, the limit of \FSe\ as $\varepsilon\rightarrow0$. This shows
that \FSe, a~variant of boosting, and the
lasso [\citet{Tibshirani:1996}] are related in some settings.
\citet{HTTW:2007} showed in general that the solution path of $\mathrm{FS}_{0}$ is
equivalent to the path of the monotone lasso.

However, much less work has focused on stagewise algorithms assuming
an arbitrary learning parameter $0<\nu\le1$. An important exception
is \citet{Buhlmann:2006} who studied \l2boost with componentwise linear
least squares, the same algorithm studied here, and
proved consistency for arbitrary $\nu$ under a sparsity assumption
where $p$ can increase at an
exponential rate relative to $n$. As pointed out in \citet
{Buhlmann:2006}, the \FSe\
algorithm studied by \citet{EfronLAR} bears similarities to \l2boost.
It is identical to Algorithm~\ref{A:L2Boost}, except for line 4,
where $\varepsilon$ is used in place of $\nu$ and
\[
\F_m = \F_{m-1} + \varepsilon\d_m\X_{k_m},\qquad
\mbox{where } \delta_m=\sign[\corr(\g_m, \X_{k_m})].
\]
Thus, \FSe\ replaces the gradient-correlation $\rho_m$ with the
sign of the gradient-correlation $\d_m$. For infinitesimally small
$\nu$ this difference appears to be inconsequential, and it is
generally believed that the two limiting solution paths are
equal [\citet{H:2007}]. In general, however, for arbitrary $0<\nu\le1$,
the two solution paths are different. Indeed, \citet{Buhlmann:2006}
indicated certain unique advantages possessed by \l2boost. Other
related work includes \citet{Buhlmann:Yu:2003}, who described a
bias-variance decomposition of the mean-squared-error of a variant of
\l2boost.

\subsection{Outline and contributions}
In this paper, we investigate the properties of \l2boost assuming an
arbitrary learning parameter $0<\nu\le1$. During \l2boost's descent
along a fixed coordinate direction, a new coordinate becomes more
favorable when it becomes closest to the current gradient. But when
does this actually occur? We provide an exact simple closed form
expression for this quantity: the number of iterations to favorability
(Theorem~\ref{criticalpoint.theorem} of Section~\ref{S:fixedDescent}).
This core identity is used to describe \l2boost's solution path
(Theorem~\ref{full.path.solution.general}), to introduce new tools for
studying its path and to study and characterize some of the
algorithm's unique properties. One of these is active set cycling, a
property where the algorithm spends lengthy periods of time cycling
between the same coordinates when~$\nu$ is small
(Section~\ref{S:cyclingBehavior}).

Our fixed descent identity also reveals how correlation affects
\l2boost's ability to select variables in highly correlated problems.
We identify a \textit{repressible condition} that prevents a new variable
from entering the active set, even though that variable may be highly
desirable (Section~\ref{S:repressibility}). Using a~data augmentation
approach, similar to that used for calculating the elastic
net [\citet{Zou:Hast:2005}], we describe a simple method for adding
$L_2$-penalization to \l2boost (Section~\ref{S:elasticBoost}). In
combination with decorrelation, this reverses the repressible
condition and improves \l2boost's performance in correlated problems.
Because \l2boost is known to approximate forward stagewise algorithms
for arbitrarily small $\nu$, it is natural to expect these results to
apply to such algorithms like LAR and lasso, and thus our results
provide a new explanation for why these algorithms may perform poorly
in correlated settings and why methods like the elastic net, which
makes use of $L_2$-penalization, are more adept in such settings. All
proofs in this manuscript can be found in the supplemental
article [\citet{EhrlingerIshwaran:2012}].

\section{Fixed descent analysis}\label{S:fixedDescent}

To analyze \l2boost we introduce the following notation useful for
describing its solution path. Let $\{l_1,\ldots,l_\Ms\}$ be the
$\Ms\le M$ nonduplicated values in order of appearance of the
selected coordinate directions $\bb_M=\{k_1,\ldots,k_M\}$. We refer to
these ordered, nonduplicated values as \textit{critical directions} of
the path. For example, if $\bb_M=\{5,5,\allowbreak 5,3,5,1,4,4,5\}$, the critical
directions are $\{5,3,5,1,4,5\}$ and $\Ms=6$. To formally describe
the solution path we introduce the following nomenclature.\looseness=-1

\begin{definition}\label{path.def}
The descent length along a critical direction $l_r$ is denoted
by~$L_r$. The critical point $S_r$ is the step number at which the
descent along~$l_r$ ends. Thus, following step $S_{r-1}$, the descent
is along $l_r$ for a total of $L_r$ steps, ending at step $S_r$.
\end{definition}

\begin{figure}

\includegraphics{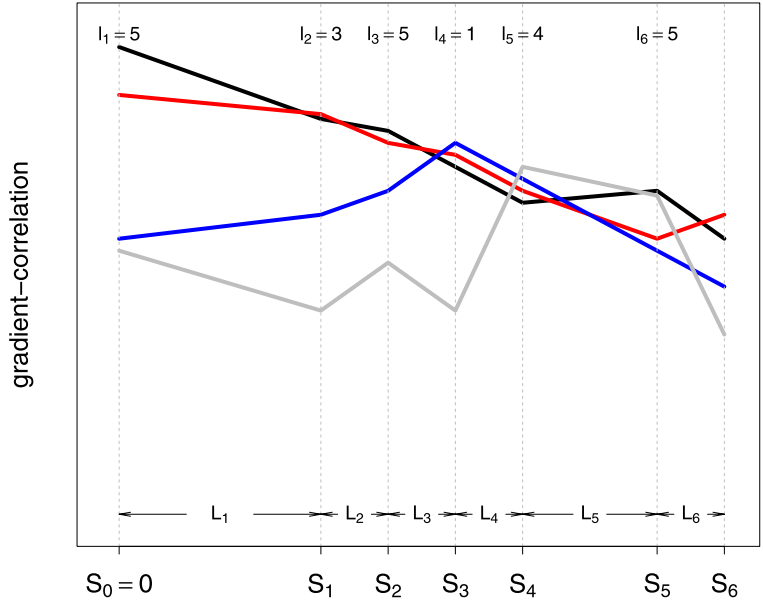}

\caption{\sl Solution path for \l2boost where $\bb_M = \{5, 5, 5, 3,
5, 1, 4, 4, 5\}$. The $M^*=6$ critical directions are $(l_r)_1^6 =
(5, 3, 5, 1, 4, 5)$ with critical descent step lengths $(L_r)_1^6 =
(3,1,1,1,2,1)$ and critical points $(S_r)_1^6 =(3,4,5,6,8,9)$.}\label{figure1}
\end{figure}

The set of values $(l_r,L_r,S_r)_{1}^{M^*}$ can be used to formally
describe the solution path of \l2boost: the algorithm
begins by descending along direction~$l_1$ (the first critical
direction) for $L_1$ steps, after which it switches to a descent along
direction $l_2$ (the second critical direction) for a total of $L_2$
steps. This continues with the last descent along $l_\Ms$ (the final
critical direction) for a total of $L_\Ms$ steps. See Figure~\ref{figure1} for
illustration of the notation.

A key observation is that \l2boost's behavior along a given descent is
deterministic except for its descent length $L_r$ (number of steps).
If we could determine the descent length, a quantity we show is highly
amenable to analysis, then an exact description of the solution path
becomes possible as \l2boost can be conceptualized as collection of
such fixed paths.

Imagine then that we are at step $m'$ of the algorithm and that in the
following step a new critical direction $k$ is formed. Let us study
the descent along $k$ for the next $m=1,\ldots,M'$ steps. Thus, in
the $m$th step of the descent along $k$, the predictor is
\[
\F_{k,m} = \F_{k,m-1} + \nu\rho_{k,m}\X_k,\qquad\mbox{where }
\rho_{k,m}=\X_k^T(\y- \F_{k,m-1}).
\]
Consider then Algorithm~\ref{A:IncrementalL2Boost} which repeatedly
boosts the predictor along the $k$th direction for a total of $M'$
steps.
\begin{algorithm}[t]
\caption{\l2boost (Fixed direction, $k$)}\label{A:IncrementalL2Boost}
\begin{algorithmic}[1]
\STATE$\F_{k,0}=\F_{m'}$
\FOR{$m=1$ to $M'$}
\STATE$\F_{k,m} = \F_{k,m-1} + \nu\rho_{k,m}\X_k$, where $\rho
_{k,m}=\X_k^T(\y- \F_{k,m-1})$
\ENDFOR
\end{algorithmic}
\end{algorithm}

The following result states a closed form solution for the $m$-step
predictor of Algorithm~\ref{A:IncrementalL2Boost} and will be crucial
to our characterization of \l2boost.

\begin{theorem}\label{incremental.operator.theorem}
$\F_{k,m} = \F_{k,0} + \nu_m\rho_{k,1}\X_k$,
where $\nu_m=1 - (1-\nu)^m$ and $\rho_{k,1}=\X_k^T(\y-\F_{k,0})$.
\end{theorem}

Theorem~\ref{incremental.operator.theorem} shows that taking a single
step with learning parameter $\nu_m$ yields the same limit as taking
$m$ steps with the smaller learning parameter~$\nu$. The result also
sheds insight into how $\nu$ slows the descent relative to stagewise
regression. Notice that the $m$-step predictor can be written as
\[
\label{slowlearning.effect}
\F_{k,m} = \underbrace{\F_{k,0} + \rho_{k,1}\X_k}_{\mathrm
{stagewise}} -
\underbrace{(1-\nu)^m\rho_{k,1}\X_k}_{\mathrm{slow\ learning}}.
\]
The first term on the right is the predictor from a greedy stagewise
step, while the second term represents the effect of slow-learning.
This latter term is what slows the descent relative to a greedy step.
When $m\rightarrow\infty$ this term vanishes, and we end up with
stagewise fitting, $\nu=1$.

\subsection{Directional change in the descent}\label{S:descentDirectionChange}
Theorem~\ref{incremental.operator.theorem} shows how to take a~large
boosting step in place of many small steps, but it does not indicate
how many steps must be taken along $k$ before a new variable enters
the solution path. If this were known, then the entire $k$-descent
could be characterized in terms of a single step.

To determine the descent length, suppose that \l2boost has descended
along $k$ for a total of $m$ steps. At step $m+1$ the algorithm must
decide whether to continue along $k$ or to select a new direction $j$.
To determine when to switch directions, we introduce the following
definition.

\begin{definition}\label{favorable.def}
A direction $j$ is said to be more favorable than $k$ at step $m+1$ if
$|\rho_{k,m}|\ge|\rho_{j,m}|$ and $|\rho_{k,m+1}|<|\rho_{j,m+1}|$.
Thus, if $j$ is more favorable at $m+1$, the descent switches to $j$
for step $m+1$.
\end{definition}

To determine when $j$ becomes more favorable, it will be useful to
have a~closed form expression for $\rho_{k,m+1}$ and $\rho_{j,m+1}$. By
Theorem~\ref{incremental.operator.theorem},
\begin{eqnarray*}\label{rho.mj.def}
\rho_{j,m+1}
&=&\X_j^T(\y-\F_{k,m})
\\
&=&\X_j^T [(\y-\F_{k,0})-\nu_{m}\rho_{k,1}\X_k ]
\\
&=&\rho_{j,1} -\nu_{m}\rho_{k,1}R_{j,k},
\end{eqnarray*}
where $R_{j,k}=\X_j^T\X_k$. Setting $j=k$ yields $\rho_{k,m+1} =
(1-\nu)^{m}\rho_{k,1}$. Therefore, $|\rho_{k,m+1}|<|\rho_{j,m+1}|$
if and only if
\[
(1-\nu)^{2m}\rho_{k,1}^2<(\rho_{j,1} -\nu_{m}\rho_{k,1}R_{j,k})^2.
\]
Dividing throughout by $\rho_{k,1}$, with a little bit of rearrangement,
this becomes
%
\begin{equation}
\label{holy.grail.ineq}
(1-\nu)^{2m} < [(1-\nu)^{m}R_{j,k} + (d_{j,k} -R_{j,k}) ]^2,
\end{equation}
where $d_{j,k}=\rho_{j,1}/\rho_{k,1}$. Notice importantly that
$|d_{j,k}|\le1$ because $k$ is the direction with maximal
gradient-correlation at the start of the descent. It is also useful
to keep in mind that $R_{j,k}$ is the sample correlation of $\X_j$ and
$\X_k$ due to \mref{data.standardization}, and thus $|R_{j,k}|\le1$.
The following result states the number of steps taken along $k$ before
$j$ becomes more favorable.

\begin{theorem}\label{criticalpoint.theorem}
The number of steps $m_{j,k}$ taken along $k$ so that $j$
becomes more favorable than $k$ at $m_{j,k}+1$ is the largest integer
$m$ such that
%
\begin{equation}
\label{critical.point}
(1-\nu)^{m-1}\ge\frac{|d_{j,k}-R_{j,k}|}{1 - R_{j,k}\sign(d_{j,k}-R_{j,k})}.
\end{equation}
It follows that for $0< \nu< 1$
%
\begin{equation}
\label{critical.step.size}
m_{j,k} =
\floor\biggl[1 + \frac{\log|d_{j,k}-R_{j,k}|
-\log(1-R_{j,k}\sign(d_{j,k}-R_{j,k}) )}{\log(1-\nu)} \biggr],\hspace*{-30pt}
\end{equation}
where $\floor(z)$ is the largest integer less than or equal to $z$.
\end{theorem}

\begin{remark}\label{repressible.remark}
In particular, notice that $m_{j,k}=\infty$ when $d_{j,k}=R_{j,k}$
[adopting the standard convention that $\sign(0)=0$ and assuming that
$\nu<1$]. We call $d_{j,k}=R_{j,k}$ the repressible condition.
Section~\ref{S:repressibility} will show that repressibility plays a
key role in \l2boost's behavior in correlated settings.
\end{remark}

\begin{remark}
When $\nu=1$ we obtain $m_{j,k}=1$
from \mref{critical.point} which corresponds to greedy stagewise
fitting. Because this makes the $\nu=1$ case uninteresting,
we shall hereafter assume that $0<\nu<1$.
\end{remark}

\subsection{Defining the solution path}

Theorem~\ref{criticalpoint.theorem} immediately shows that the problem
of determining the next variable to enter the solution path can be
recast as finding the direction requiring the fewest number of steps
$m_{j,k}$ to favorability. When combined with
Theorem~\ref{incremental.operator.theorem}, this characterizes the
entire descent and can be used to characterize \l2boost's solution
path.\vadjust{\goodbreak}

As before, assume that $k$ corresponds to the first critical direction
of the path, that is, $l_1=k$. By Theorem~\ref{criticalpoint.theorem},
\l2boost descends along $k$ for a total of $S_1=L_1$ steps, where
$L_1=m_{l_2,k}$ and $l_2$ is the coordinate requiring the smallest
number of steps to become more favorable than $k$. By
Theorem~\ref{incremental.operator.theorem}, the predictor at step
$S_1$ is
\[
\F_{S_1} = \F_0 + \nu_{L_1}\rho_{l_1}^{(1)}\X_{l_1},\qquad
\mbox{where }\rho_{l_1}^{(1)}=\X_{l_1}^T(\y-\F_0).
\]
Applying Theorem~\ref{incremental.operator.theorem} once again, but
now using a descent along $l_2$ initialized at $\F_{S_1}$, and
continuing this argument recursively, as well as using the
representation for the number of steps from
Theorem~\ref{criticalpoint.theorem}, yields
Theorem~\ref{full.path.solution.general}, which presents a recursive
description of \l2boost's solution path.

\begin{theorem}\label{full.path.solution.general}
$\F_{S_r} = \F_{S_{r-1}} + \nu_{L_r}\rho_{l_r}^{(r)}\X_{l_r}$,
where $\{(l_r,L_r,S_r,\rho_{l_r}^{(r)})\}_1^{\Ms}$ are
determined recursively from
\begin{eqnarray*}
l_1&=&\mathop\argmax_{1\le j\le p} |\X_j^T(\y-\F_0)|,\qquad
l_{r+1}=\mathop\argmax_{j\ne l_r}\bigl|\rho_j^{(r+1)}\bigr|,
\\
M_j^{(r)} &=&
\floor\biggl[1 + \frac{ \log|D_j^{(r)} - R_{j,l_r}|
- \log(1 - R_{j,l_r}\sign(D_j^{(r)} - R_{j,l_r})
)}{\log(1-\nu)} \biggr],
\\
L_r&=&M_{l_{r+1}}^{(r)},\qquad S_r = S_{r-1} + L_r,\qquad S_0=0,
\\
D_j^{(r)} &=& \frac{\rho_j^{(r)}}{\rho_{l_r}^{(r)}},\qquad
\rho_{j}^{(r+1)} = \X_j^T(\y-\F_{S_{r}}) = \rho_{j}^{(r)}
- \nu_{L_r}\rho_{l_r}^{(r)}R_{j,l_r}.
\end{eqnarray*}
\end{theorem}

\begin{remark}\label{step.number.tie.remark}
A technical issue arises in Theorem~\ref{full.path.solution.general}
when $M_j^{(r)}$ is not\vspace*{2pt} unique. Non-uniqueness can occur due to
rounding which is caused by the floor function used in the definition
of $m_{j,k}$. This is why line 1 selects the next critical
value, $l_{r+1}$, by maximizing the absolute gradient-correlation
$|\rho_j^{(r+1)}|$ and not by minimizing the step number $M_j^{(r)}$.
This definition for~$l_{r+1}$ is equivalent to the two-step solution
\[
l_{r+1}\leftarrow\mathop\argmax_{j\in
l_{r+1}}\bigl|\rho_j^{(r+1)}\bigr|,\qquad
\mbox{ where }
l_{r+1}=\mathop\argmin_{j\ne l_r} \bigl\{M_j^{(r)}\bigr\}.
\]
\end{remark}

\begin{remark}\label{equal.gradient.remark} Another technical issue arises
when there is a tie in the absolute gradient-correlation. In line 3
of Algorithm~\ref{A:L2Boost} it may be possible for two coordinates,
say $j$ and $k$, to have equal gradient-correlations at step $m>1$.
Theorem~\ref{full.path.solution.general} implicitly deals with such
ties due to Definition~\ref{favorable.def}. For example, suppose that
the first $m-1$ steps are along $k$ with the tie occurring at step $m$.
In the language of Theorem~\ref{criticalpoint.theorem}, because $j$
becomes more favorable than $k$ at $m+1$, where $m=m_{j,k}$, we have
\[
|\rho_{j,m-1}| < |\rho_{k,m-1}|,\qquad
|\rho_{j,m}| = |\rho_{k,m}|,\qquad
|\rho_{j,m+1}| > |\rho_{k,m+1}|.\vadjust{\goodbreak}
\]
In this example, Theorem~\ref{full.path.solution.general} resolves
the tie at $m$ by continuing to descend along $k$, then switching to
$j$ at step $m+1$. Although Algorithm~\ref{A:L2Boost} does not explicitly
address this issue, the potential discrepancy is minor because such
ties should rarely occur in practice. This is because for
$|\rho_{j,m}|=|\rho_{k,m}|$ to hold, the value inside the floor
function of \mref{critical.step.size} used to define $m_{j,k}$ must be
an integer (a careful analysis of the proof of
Theorem~\ref{criticalpoint.theorem} shows why). A tie can occur only
when this value is an integer which is numerically unlikely to
occur.\looseness=-1\vspace*{-3pt}
\end{remark}

\begin{remark}
Theorem~\ref{full.path.solution.general} immediately yields
a recursive solution for the coefficient vector, $\b$.
The solution path for $\b$ is the piecewise solution
\[
\b^{(r)} = \b^{(r-1)} + \nu_{L_r}\rho_{l_r}^{(r)}\one_{l_r},\qquad
\b^{(0)}=\zero,
\]
where $\one_{l_r}\in\RR^p$ is the vector with one in coordinate $l_r$
and zero elsewhere.\vspace*{-3pt}
\end{remark}

\subsection{Illustration: Diabetes data}

Aside from the technical issue of ties,
Theorem~\ref{full.path.solution.general} and Algorithm~\ref{A:L2Boost}
are equivalent. For convenience, we state
Theorem~\ref{full.path.solution.general} in an algorithmic form to
facilitate comparison with Algorithm~\ref{A:L2Boost}; see
Algorithm~\ref{A:L2BoostPath}. Computationally,
Algorithm~\ref{A:L2BoostPath} improves upon Algorithm~\ref{A:L2Boost}
by avoiding taking many small steps along a given descent. However,
the difference is not substantial because the benefits only apply when
$\nu$ is small, and as we will show later (Section~\ref{S:cyclingBehavior}), this forces the
algorithm to cycle between its variables following the first descent,
thus mitigating its ability to take large steps. Thus, strictly
speaking, the benefit of Algorithm~\ref{A:L2BoostPath} is confined
primarily to the first descent.

\begin{algorithm}[b]
\caption{\l2boost (Solution path)}\label{A:L2BoostPath}
\begin{algorithmic}[1]
\STATE$\F_0=\zero$; $S_0=0$; $l_1=\argmax_{1\le j\le p} |\X_j^T\y|$
\FOR{$r=1$ to $\Ms$}
\STATE$l_{r+1}=\argmax_{j\ne l_r}|\rho_j^{(r+1)}|$;
$\rho_{j}^{(r+1)} = \rho_{j}^{(r)} - \nu_{L_r}\rho_{l_r}^{(r)}R_{j,l_r}$
\STATE$L_r=M_{l_{r+1}}^{(r)}$; $S_r=S_{r-1}+L_r$

\STATE$\F_{S_r} = \F_{S_{r-1}} + \nu_{L_r}\rho_{l_r}^{(r)}\X_{l_r}$
\ENDFOR
\end{algorithmic}
\end{algorithm}

To investigate the differences between the two algorithms
we analyzed the diabetes data used in \citet{EfronLAR}. The data
consists of $n=442$ patients in which the response of interest, $y$,
is a quantitative measure of disease progression for a patient. In
total there are 64 variables, that includes 10 baseline measurements
for each patient, 45 interactions and 9 quadratic terms.

\begin{figure}

\includegraphics{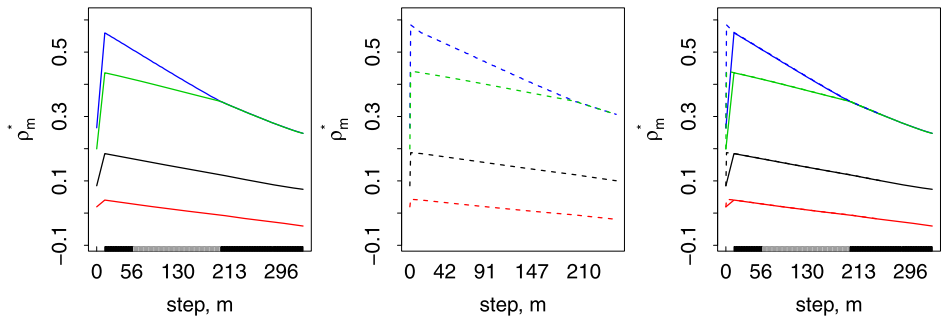}

\caption{\l2boost applied to the diabetes data. First two
panels display standardized gradient\-correlation $\rho_m^*$
against step number $m$ for Algorithms~\protect\ref{A:L2BoostPath}
and~\protect\ref{A:L2Boost}, respectively. Only coordinates in the solution
path are displayed (a total of four). The third panel superimposes the
first two
panels. All analyses used $\nu=0.005$.}\label{figure2}\vspace*{-3pt}
\end{figure}

In order to compare results, we translated each iteration, $r$, used
by Algorithm~\ref{A:L2BoostPath} into its corresponding number of
steps, $m$. Thus, while we ran Algorithm~\ref{A:L2BoostPath} for
$M^*=250$ iterations, this translated into $M=332$ steps. As
expected, this difference is primarily due to the first iteration
$r=1$ which took\vadjust{\goodbreak} $m=14$ steps along the first critical direction
(first panel of Figure~\ref{figure2}; the rug indicates critical
points, $S_r$).
There are other instances where Algorithm~\ref{A:L2BoostPath} took
more than one step (corresponding to the light grey tick marks on the
rug), but these were generally steps of length 2.
The standardized gradient-correlation is plotted along the $y$-axis of
the figure.
The standardized
gradient-correlation for step $m$ was defined as (using the notation
of Algorithm~\ref{A:L2Boost})
%
\begin{equation}
\label{stand.grad.corr}
\rho_{m}^*
=\frac{\X_{k_m}^T\g_m}{\sqrt{\X_{k_m}^T\X_{k_m}}
\sqrt{\g_{m}^T\g_{m}}} =
\frac{\rho_m}{\sqrt{\g_{m}^T\g_{m}}}.
\end{equation}
The middle panel displays the results using Algorithm~\ref{A:L2Boost}
with $M=250$ steps. Clearly, the greatest gains from
Algorithm~\ref{A:L2BoostPath} occur along the $r=1$ descent. One can
see this most clearly from the last panel which superimposes the first
two panels.\vspace*{-3pt}

\begin{remark}
Note a potential computational optimization exists in
Algorithm~\ref{A:L2BoostPath}. It is possible to calculate the
correlation values only once as each new variable enters the active
set, then cache these values for future calculations. Thus, when
$l_{r+1}$ is a new variable in the active set, we calculate
$(R_{j,l_{r+1}})_{j=1}^p$. The updated gradient-correlation is
calculated efficiently by using addition and scalar multiplication
using the previous gradient-correlation and the cached correlation
coefficients
\[
\rho_{j}^{(r+1)} = \rho_{j}^{(r)} -
\nu_{L_r}\rho_{l_r}^{(r)}R_{j,l_r}.
\]
This is in contrast to
Algorithm~\ref{A:L2Boost} which requires a vector multiplication of
dimension $p$ at each step $m$ to update the gradient-correlation:
$\rho_m = \X_{k_m}^T\g_m$.\vspace*{-3pt}
\end{remark}

\begin{remark}\label{active.set.remark}
Above, when we refer to the ``active set,'' we mean the unique set of
critical directions in the current solution path. This term will be used
repeatedly throughout the paper.\vspace*{-3pt}
\end{remark}

\subsection{Visualizing the solution path}
Throughout the paper we illustrate different ways of utilizing
$m_{j,k}$ of Theorem~\ref{criticalpoint.theorem} to explore \l2boost.
So far we have confined the use of Theorem~\ref{criticalpoint.theorem}
to determining the descent length along a fixed direction, but another
interesting application is determining how far a given variable is
from the active set. Note that although
Theorem~\ref{criticalpoint.theorem} was described in terms of an
active set of only one coordinate, it applies in general, regardless of
the size of the active set. Thus, $m_{j,k}$ can be calculated at any
step $m$ to determine the number of steps required for $j$ to become
more favorable than the current direction, $k$. This value represents
the distance of $j$ to the solution path and can be used to visualize
it.

\begin{figure}

\includegraphics{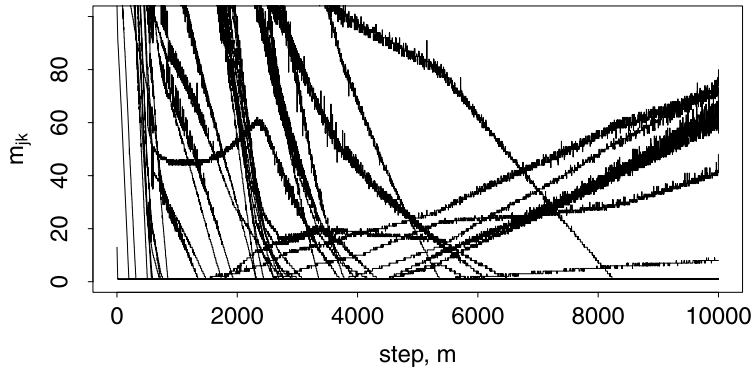}

\caption{Distance $m_{j,k}$ of each variable $j$ to favorability
relative to the current descent~$k$ (results based on
Algorithm~\protect\ref{A:L2Boost} where $\nu=0.005$). For visual clarity
the $m_{j,k}$ values have been smoothed using a running median
smoother.}\label{figure3}
\end{figure}

To demonstrate this, we applied Algorithm~\ref{A:L2Boost} to the
diabetes data for $M=10\mbox{,}000$ steps and recorded $m_{j,k}$ for each of
the $p=64$ variables. Figure~\ref{figure3} records these values.
Each ``jagged path'' in the figure is the trace over the 10,000 steps
for a variable $j$. Each point on the path equals the number of steps
$m_{j,k}$ to favorability relative to the current descent $k\ne j$.
The patterns are quite interesting. The top variables have $m_{j,k}$
values which quickly drop within the first 1000 steps. Another group
of variables have values which take much longer to drop, doing so
somewhere between 2000 to 4000 steps, but then increase almost
immediately. These variables enter the solution path but then quickly
become unattractive regardless of the descent direction.

It has become popular to visualize the solution path of forward
stagewise algorithms by plotting their gradient-correlation paths
and/or their coefficient paths. Figure~\ref{figure3} is a similar
tool. A unique
feature of $m_{j,k}$ is that it depends not only on the
gradient-correlation (via $d_{j,k}$), but also the
correlation in the $x$-variables (via $R_{j,k}$) and the learning
parameter $\nu$. In this manner, Figure~\ref{figure3} offers a new
tool for
understanding and exploring such algorithms.

\section{Cycling behavior}\label{S:cyclingBehavior}

It has been widely observed that decreasing the regularization
parameter slows the convergence of stagewise descent
algorithms. \citet{EfronLAR} showed that the \FSe\ algorithm tracks
the equiangular direction of the LAR path for arbitrarily small
$\varepsilon$.
To achieve what LAR does in a single step, the \FSe\ algorithm may
require thousands of small steps in a direction tightly clustered
around the equiangular vector, eventually ending up at nearly the same
point as LAR.

We show that \l2boost exhibits this same phenomenon. We do so by
describing this property as an active set cycling phenomenon. Using
results from the earlier fixed descent analysis, we show in the case
of an active set of two variables that \l2boost systematically
switches (cycles) between its two variables when $\nu$ is small. For
an arbitrarily small $\nu$ this forces the absolute
gradient-correlations for the active set variables to be nearly equal.
This point of equality represents a~singularity point that triggers a~near-perpetual deterministic cycle between the variables, ending only
when a~new variable enters the active set with nearly the same
absolute gradient-correlation.

\subsection{\l2boost's gradient equality point}
Our insight will come from looking at
Theorem~\ref{criticalpoint.theorem} in more depth. As before, assume
the algorithm has been initialized so that $k$ is the first critical
step. Previously the descent along~$k$ was described in terms of
steps, but this can be equivalently expressed in units of the ``step
size'' taken. Define
\[
\nu_{j,k}=\nu_{m_{j,k}}=1-(1-\nu)^{m_{j,k}}.
\]
Recall that Theorem~\ref{incremental.operator.theorem} showed that
a single step along $k$ with $\nu$ replaced with $\nu_{j,k}$
yields the same limit as $m_{j,k}$ steps along $k$ using $\nu$. We
call $\nu_{j,k}$ the step size taken along $k$. Because $j$ becomes
more favorable than $k$ at $m_{j,k}+1$, the gradient following a step
size of $\nu_{j,k}$ along $k$ satisfies
%
\begin{equation}
\label{overshoot}
|\X_j^T (\y-\F_0-\nu_{j,k}\rho_{k,1}\X_k ) |
<
|\X_k^T (\y-\F_0-\nu_{j,k}\rho_{k,1}\X_k ) |.
\end{equation}
This applies to all coordinates $j\ne k$, and in particular holds for
the second critical direction, $l_2$, which rephrased in terms of
step size, is the smallest~$\nu_{j,k}$ value,
\[
l_2=\mathop\argmin_{j\ne k} \{\nu_{j,k}\}.
\]
Although inequality \mref{overshoot} is strict, it becomes
arbitrarily close to equality with shrinking $\nu$. With a little bit
of rearranging, \mref{critical.point} implies that
%
\begin{equation}
\label{nuhat.solution}
\nuhat_j< \nu_{j,k},\qquad\mbox{where }
\nuhat_j
= 1 - \frac{|d_{j,k} - R_{j,k}|}
{1 - R_{j,k}\sign(d_{j,k} - R_{j,k})}.
\end{equation}
We will show $\nuhat_j$ is the step size making the
absolute gradient-correlation between $j$ and $k$ equal
%
\begin{equation}
\label{singularity.point}
|\X_j^T (\y-\F_0-\nuhat_j\rho_{k,1}\X_k ) |
=
|\X_k^T (\y-\F_0-\nuhat_j\rho_{k,1}\X_k ) |.\vadjust{\goodbreak}
\end{equation}
The next theorem shows that $\nu_{l_2,k}$ converges to the
smallest $\nuhat_j$ satisfying~\mref{singularity.point};
thus, \mref{overshoot} becomes an equality in the limit.
For convenience, we define $\nu_{j,k}^-=\nu_{m_{j,k}-1}$.\vspace*{-3pt}

\begin{theorem}\label{dynamic.nu.size}
Let $\rhohat_j = \X_j^T(\y-\F_0-\nuhat_j\rho_{k,1}\X_k)$. Then
$|\rhohat_j|=|\rhohat_k|$. Furthermore, if $ \lstar=\argmin_{j\ne k}
\{\nuhat_j\}$ and $\nuhat=\nuhat_{\lstar}$, then
$\nu_{l_2,k}^-\le\nuhat<\nu_{l_2,k}$ and
$\nu_{l_2,k}\rightarrow\nuhat$ as $\nu\rightarrow0$.\vspace*{-3pt}
\end{theorem}

Therefore, for arbitrarily small $\nu$, $\nu_{l_2,k}\asymp\nuhat$ and
$k$ and $l_2$ will have near-equal absolute gradient-correlations.
This latter property triggers two-cycling. To see why, let us assume
for the moment that the active set variables have equal absolute
gradient-correlations. Then by a direct application of
Theorem~\ref{criticalpoint.theorem}, one can show that the number of
steps taken along $l_2$ before $k$ becomes more favorable is $m=1$.
Thus, following the descent along $k$, the algorithm switches to
$l_2$, but then immediately switches back to $k$. If $\nu$ is small
enough, this process is repeated, setting off a two-cycling pattern.

The next result is a formal statement of these arguments.
Define
\[
d_{j,k}^{(m)}=\frac{\rho_{j,m}}{\rho_{k,m}},\qquad\mbox{where }
\rho_{l,m}=\X_l^T(\y-\F_{m-1}), 1\le l \le p.
\]
For notational convenience, let $j=l_2$ and $m=m_{j,k}$. For
technical reasons we shall assume $d_{j,k}^{(m)}\ne R_{j,k}$.\vspace*{1pt} Recall
Remark~\ref{repressible.remark} showed that $d_{j,k}^{(m)}= R_{j,k}$,
the repressible condition, yields an infinite number of steps to
favorability. Thus, for $k$ to be even eligible for favorability
we must have $d_{j,k}^{(m)}\ne R_{j,k}$.\vspace*{-3pt}

\begin{theorem}\label{long.descent.followed.two.cycle}
If the first two critical directions are $(k,j)$ and
$\nu_{j,k}=\nuhat_j$, then $k$ is favored over $j$ for the next step
after $j$ if $d_{j,k}^{(m)}\ne R_{j,k}$.\vspace*{-3pt}
\end{theorem}

Theorem~\ref{long.descent.followed.two.cycle} assumes that
$\nu_{j,k}=\nuhat_j$. While this only holds in the limit, the two
values should be nearly equal for arbitrarily small $\nu$, and thus
the assumption is reasonable. Notice also that
Theorem~\ref{long.descent.followed.two.cycle} only shows that $k$ is
more favorable than~$j$, and not that the algorithm switches to $k$.
However, we can see that this must be the case. For arbitrarily small
$\nu$, $k$'s gradient-correlation should be nearly equal to $j$'s, and by
definition, $j$ has maximal absolute gradient-correlation along the
second descent.

Indeed, the following result shows that the absolute
gradient-correlations for $k$ and $j$ can be made arbitrarily close
for small enough $\nu$ for any step $r\ge1$ following the descent
along $k$. The result also shows that the sign of the
gradient-correlation is preserved when $\nu$ is arbitrarily small, a
fact that we shall use later.\vspace*{-3pt}

\begin{theorem}\label{gradient.correlation.equality}
$\rho_{j,m+r}/\rho_{k,m+r}\rightarrow\sign(\rhohat_j)/\sign
(\rhohat_k)$
as $\nu\rightarrow0$ for each $r\ge1$.\vspace*{-3pt}
\end{theorem}

Combining Theorems~\ref{long.descent.followed.two.cycle}
and~\ref{gradient.correlation.equality}, we see that if $\nu$ is
small enough, the first three critical directions of the\vadjust{\goodbreak} path must
be $(k,j,k)$ with critical points $(m,m+1, m+2)$. And once the descent
switches back to $k$, it is clear from the same argument that
the next critical direction, $l_4$, will be $j$, and so forth.

\subsection{Illustration of two-cycling}

We present a numerical example demonstrating two-cycling. For our
example, we simulated data according to
\[
\y= \X\b+ \eb,\qquad \eb\sim N(\zero,\I),
\]
where $n = 100$, and $p = 40$. The first 10 coordinates of $\b$ were
set to 5, with the remaining coordinates set to 0. The design matrix
$\X$ was simulated by drawing its entries independently from a standard
normal distribution.

\begin{figure}[b]

\includegraphics{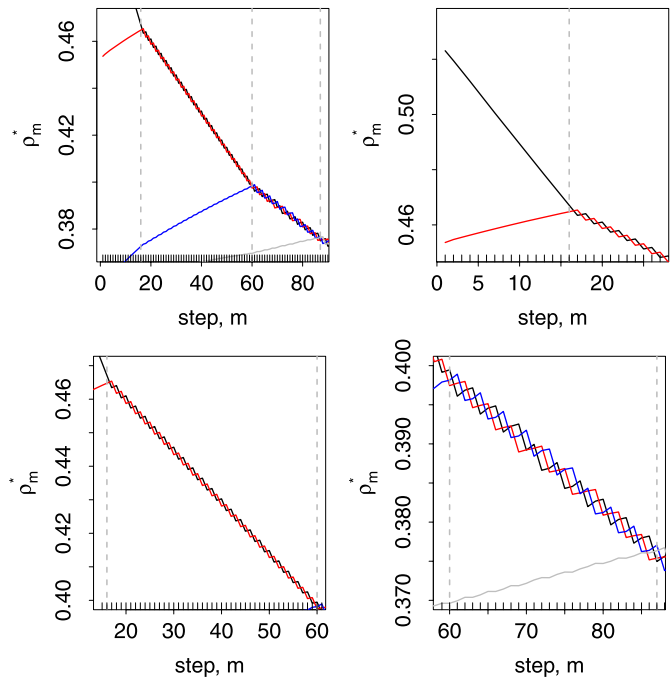}

\caption{Standardized gradient-correlation path for $\nu=0.01$. Top
left panel details the path through the first three active variables,
the remaining panels detail each active variable descent.}\label{figure4}
\end{figure}

Figure~\ref{figure4} plots the standardized
gradient-correlations \mref{stand.grad.corr} from
Algorithm~\ref{A:L2BoostPath} using $\nu=0.01$. As done earlier, we
have converted iterations $r$ into step numbers $m$ along the
$x$-axis. The plots show the behavior of each coordinate within an
active set descent. The rug marks show each step $m$ for clarity, and
dashed vertical lines indicate the step $m_{j,k}$ where the next step
adds a new critical direction to the solution path. The top left
panel shows the complete descent along the first three active
variables. The remaining panels detail the coordinate behavior as the
active set increases from one to three coordinates.

The top right panel shows repeated selection of the $l_1$ direction
shown in black. The last step along $l_1$ occurs at $m_{j,k}$ marked
with the vertical dashed line, where the next step is along the $l_2$
direction shown in red. This point marks the beginning of the
two-cycling behavior, which continues in the lower left panel. At each
step, the algorithm systematically switches between the
$l_1$ and $l_2$ directions, until an additional direction becomes more
favorable. The cycling pattern is $\{l_1,l_2,l_1,l_2,\ldots\}$. The
lower right panel demonstrates three-cycling behavior. Here it is
instructive to note that the order of selection within three-cycling
is nondeterministic. In this panel the order starts as $\{l_3, l_2,
l_1, \ldots\}$, but changes near $m=70$ to $\{\ldots, l_3, l_1,l_2,
\ldots\}$. As discussed later, nondeterministic cycling patterns are
typical behavior of higher order cycling (active sets of size greater
than two).

\subsection{The limiting path}

Here we provide a formal limiting result of two-cycling. The result
can be viewed as the analog of Theorem~\ref{dynamic.nu.size} when the
active set involves two variables. Using a slightly modified version
of \l2boost we show that for arbitrarily small $\nu$, if the algorithm
cycles between its two active variables, it does so until a new
variable enters the active set with the same absolute
gradient-correlation.

Assume the active set is $\aa=\{k,j\}$ and that $k$ and
$j$ are cycling according to $(k,j,k,j,\ldots)$. The $m$-step
predictor for $m=1,\ldots,M$ is
%
\begin{equation}\label{kj.cycling}
\F_{m} =
\cases{
\F_{m-1} + \nu\rho_{k,m}\X_k,
&\quad if $m$ is odd,
\cr
\F_{m-1} + \nu\rho_{j,m}\X_j,
&\quad if $m$ is even,}
\end{equation}
where $\rho_{l,m}=\X_l^T(\y-\F_{m-1})$. The cycling
pattern \mref{kj.cycling} is assumed to persist for a minimum length
of $M\ge3$.

It will simplify matters if the cycling is assumed to be initialized
with strict equality of the gradient correlations:
$|\rho_{k,1}|=|\rho_{j,1}|$. With an arbitrarily small $\nu$, this
will force near equal absolute gradient-correlations at each step and
by Theorem~\ref{gradient.correlation.equality} will preserve the sign
of the gradient-correlation. We assume
\[
\frac{ \rho_{j,m}}{\rho_{k,m}} =
\frac{\sign(\rho_{j,1})}{\sign(\rho_{k,1})}\qquad
\mbox{for } m\ge1.
\]
It should be emphasized that the above assumptions represent a
simplified version of \l2boost. In practice, we would have
\[
\rho_{j,m}=s\rho_{k,m} + O(\nu),
\]
where $s=\sign(\rho_{j,1})/\sign(\rho_{k,1})$. However, for
convenience we will not concern ourselves with this level of detail
here. Readers can consult \citet{Ehrlinger:2011} for a more refined
analysis.

One way to ensure $|\rho_{k,1}|=|\rho_{j,1}|$ is to initialize the
algorithm with the limiting predictor $\F_0+\nuhat_j\rho_{k,1}\X_k$ of
Theorem~\ref{dynamic.nu.size} obtained by letting \mbox{$\nu\rightarrow0$}
along the $k$-descent. With a slight abuse of notation denote this
initial estimator by $\F_0$.\vadjust{\goodbreak} However, the fact that this specific
$\F_0$ is used does not play a direct role in the results. Under the
above assumptions, the following closed form expression for the
$m$-step predictor under two-cycling holds.

\begin{theorem}\label{two.cycle.predictor}
Assume that $\rho_{j,m}=s\rho_{k,m}$ for $m\ge1$. If $d_{j,k}\ne
R_{j,k}$, then for any $0<\nu<1/2$ satisfying $1+sR_{j,k}>\nu R_{j,k}^2$,
we have for each $m\ge1$,
\[
\F_{m} =
\cases{
\F_{0} + V_{m+1}\rho_{k,1}\biggl[\X_k
+\displaystyle\frac{V_{m-1}}{V_{m+1}}(s-\nu R_{j,k})\X_j\biggr],
&\quad if $m$ is odd,\vspace*{2pt}
\cr
\F_{0} + V_{m}\rho_{k,1} [\X_k+(s-\nu R_{j,k})\X_j ],
&\quad if $m$ is even,}
\]
where $V_m=\nu\nu_{\aa}^{-1}[1-(1-\nu_{\aa})^{m/2}]$
and $\nu_{\aa} = \nu(1 + s R_{j,k} - \nu R_{j,k}^2)$.
Note that $0<\nu_{\aa}<1$ under the asserted conditions.
\end{theorem}

To determine the above limit requires first determining when a
new direction $l\notin\aa$ becomes more favorable. For $l$
to be more favorable at $m+1$, we must have
$|\rho_{j,m+1}|<|\rho_{l,m+1}|$ when $m$ is odd, or
$|\rho_{k,m+1}|<|\rho_{l,m+1}|$ when $m$ is even. The following
result determines the number of steps to favorability. For simplicity
only the case when $m$ is odd is considered, but this does not affect
the limiting result.

\begin{theorem}\label{criticalpoint.twocycle.theorem} Assume the same
conditions as Theorem~\ref{two.cycle.predictor}.
Then $l$ becomes
more favorable than $j$ at step $m+1$ where $m$
is the largest odd integer $m\ge3$ such that
%
\begin{equation}
\label{odd.critical.point}
(1-\nu_{\aa})^{(m-1)/2} \ge
\frac{|d_{l,k}-R_{j,k,l}|}{1 - R_{j,k,l}\sign(d_{l,k} - R_{j,k,l})}
\end{equation}
where $d_{l,k}=\rho_{l,1}/\rho_{k,1}$ and
\[
R_{j,k,l} =
\frac{R_{l,k}+(s-\nu R_{j,k})R_{l,j}}{1 + s R_{j,k} - \nu R_{j,k}^2}.
\]
\end{theorem}

Clearly \mref{odd.critical.point} shares common features
with \mref{critical.point}. This is no coincidence. The bounds are
similar in nature because both are derived by seeking the
point where the absolute gradient-correlation between sets of
variables are equal. In the case of two-cycling, this is the
singularity point where $k$, $j$ and~$l$ are all equivalent in terms of
absolute gradient-correlation. The following result states the
limit of the predictor under two-cycling.

\begin{theorem}\label{dynamic step.two.cycles}
Under the conditions of Theorem~\ref{two.cycle.predictor}, the limit
of $\F_m$ as $\nu\rightarrow0$ at the next critical direction
$\lstar$
equals
\[
\Fhat= \F_0 + \nuhat\rho_{k,1} [\X_k +
s\X_j ],
\]
where $ \lstar=\argmin_{l\notin\aa}\{\nuhat_l\}$,
$\nuhat=\nuhat_{\lstar}$,
%
\begin{equation}
\label{holy.grail.twostepsize}
\nuhat_l = \biggl( 1 -
\frac{|d_{l,k}-\Rhat_{j,k,l}|}
{1 - \Rhat_{j,k,l}\sign(d_{l,k} - \Rhat_{j,k,l})}
\biggr)(1 + s R_{j,k} )^{-1},\vadjust{\goodbreak}
\end{equation}
and $\Rhat_{j,k,l} =
(R_{l,k}+ s R_{l,j})/(1 + s R_{j,k})$. Furthermore,
$|\rhohat_{\lstar}|=|\rhohat_k|=|\rhohat_j|$,
where for each $l$, $\rhohat_l = \X_l^T(\y-\Fhat)$.
\end{theorem}

This shows that the predictor moves along the combined
direction $\X_k + s\X_j$ taking a step size $\nuhat$ that makes the
absolute gradient-correlation for $\lstar$ equal to that of the active
set $\aa=\{k,j\}$. Theorem~\ref{dynamic step.two.cycles} is a direct
analog of Theorem~\ref{dynamic.nu.size} to two-cycling.

Not surprisingly, one can easily show that this limit coincides with
the LAR solution. To show this, we rewrite $\Fhat$ in a form
comparable to LAR,
\[
\Fhat= \F_0 + \nuhat|\rho_{k,1}| [\sign(\rho_{k,1})\X_k +
\sign(\rho_{j,1})\X_j ].
\]
Recall that LAR moves the shortest distance along the equiangular
vector defined by the current active set until a new variable with
equal absolute gradient-correlation is reached. The term in square
brackets above is proportional to this equiangular vector. Thus,
since $\Fhat$ is obtained by moving the shortest distance along the
equiangular vector such that $\{j,k,\lstar\}$ have equal absolute
gradient-correlation, $\Fhat$ must be identical to the LAR solution.

\subsection{General cycling}
Analysis of cycling in the general case where the active set
$\aa=\{k_i\}_{i=1}^d$ is comprised of $d\ge2$ variables is more
complex. In two-cycling we observed cycling patterns of the form
$(l_1,l_2,l_1,l_2,\ldots)$, but when $d>2$, \l2boost's cycling
patterns are often observed to be nondeterministic with no discernible
pattern in the order of selected critical directions. Moreover, one
often observes some coordinates being selected more frequently than
others.

A study of $d$-cycling has been given by \citet{Ehrlinger:2011}.
However, the analysis assumes deterministic cycling of the form
\[
(l_1,l_2,\ldots, l_d,l_{d+1}, \ldots) =
(k_1, k_2, \ldots, k_d, k_1, \ldots),
\]
which is the natural extension of the two-cycling just studied. To
accommodate this framework, a modified \l2boost procedure involving
coordinate-dependent step sizes was used. This models \l2boost's
cycling tendency of selecting some coordinates more frequently by
using the size of a step to dictate the relative frequency of
selection. Under constraints to the coordinate step sizes, equivalent
to solving a system of linear equations defining the
equiangular vector used by LAR, it was shown that the modified
\l2boost procedure yields the LAR solution in the limit. Interested
readers should consult \citet{Ehrlinger:2011} for details.

\section{Repressibility affects variable selection in correlated
settings}\label{S:repressibility}

Now we turn our attention to the issue of correlation. We have shown
that regardless of the size of the
active set a new direction $j$ becomes more favorable than the current\vadjust{\goodbreak}
direction $k$ at step $m_{j,k}+1$ where $m_{j,k}$ is the smallest
integer value satisfying
%
\begin{equation}
\label{heuristic.repressibility}
1 - \frac{|d_{j,k} - R_{j,k}|}
{1 - R_{j,k}\sign(d_{j,k} - R_{j,k})}
< 1-(1-\nu)^{m_{j,k}}.
\end{equation}
Using our previous notation, let $\nuhat_j$ and $\nu_{j,k}$ denote
the left and right-hand sides of the above inequality, respectively.

Generally, large values of $m_{j,k}$ are designed to hinder
noninformative variables from entering the solution path. If $j$
requires a large number of steps to become favorable, it is
noninformative relative to the current gradient and therefore
unattractive as a candidate. Surprisingly, however, such an
interpretation does not always apply in correlated problems. There
are situations where $j$ is informative, but $m_{j,k}$ can be
artificially large due to correlation.

To see why, suppose that $j$ is an informative variable with a
relatively large value of $d_{j,k}$. Now, if $j$ and $k$ are
correlated, so much so that $R_{j,k}\approx d_{j,k}$, then
$|d_{j,k}-R_{j,k}|\approx0$. Hence, $m_{j,k}\approx\infty$ and
$\nu_{j,k}\approx1$ due to \mref{heuristic.repressibility}. Thus,
even though $j$ is promising with a large gradient-correlation, it is
unlikely to be selected because of its high correlation with $k$.

The problem is that $j$ becomes an unlikely candidate for selection
when~$d_{j,k}$ is close to $R_{j,k}$. In fact, $m_{j,k}=\infty$ when
$d_{j,k}=R_{j,k}$ so that $j$ can never become more favorable than $k$
when the two values are equal. We have already discussed the
condition $d_{j,k}=R_{j,k}$ several times now, and have referred to it
as the \textit{repressible condition}. Repressibility plays an important
role in correlated settings. We distinguish between two types of
repressibility: weak and strong repressibility. Weak repressibility
occurs in the trivial case when $|R_{j,k}|=1$. Weak repressibility
implies that $|d_{j,k}|=|R_{j,k}|=1$. Hence the gradient-correlation
for~$j$ and~$k$ are equal in absolute value and $j$, and~$k$ are
perfectly correlated. This trivial case simply reflects a numerical
issue arising from the redundancy of the $j$ and $k$ columns of the
$\X$ design matrix. The stronger notion of repressibility, which we
refer to as strong repressibility, is required to address the
nontrivial case $|R_{j,k}|\ne1$ in which $j$ is repressed without
being perfectly correlated with $k$. The following definition
summarizes these ideas.

\begin{definition}\label{repressible.def}
We say $j$ has the strong
repressible condition if $d_{j,k}=R_{j,k}$ and $|R_{j,k}|<1$. We say
that $j$ is (strongly) repressed by $k$ when this happens. On the
other hand, $j$ has the weak repressible condition if $j$ and~$k$ are
perfectly correlated ($|R_{j,k}|=1$) and $d_{j,k}=R_{j,k}$.
\end{definition}

\subsection{An illustrative example}
We present a numerical example of how repressibility can hinder
variables from being selected. For our illustration we use example
(d) of Section 5 from \citet{Zou:Hast:2005}. The data was simulated
according to
\[
\y=\X\b+ \s\eb,\qquad \eb\sim N(\zero,\I),
\]
where $n=100$, $p=40$ and $\s=15$. The first 15 coordinates of $\b$
were set to~3; all other coordinates were 0. The design
matrix $\X=[\X_1,\ldots,\X_{40}]_{100\times40}$ was simulated
according to
%
\begin{eqnarray}
\label{hastie.sim}
\X_j & = & \Z_1 + \tau\eb_j,\qquad j = 1,\ldots,5,\nonumber
\\
\X_j & = & \Z_2 + \tau\eb_j,\qquad j = 6,\ldots,10,\nonumber
\\[-8pt]
\\[-8pt]
\X_j & = & \Z_3 + \tau\eb_j,\qquad j = 11,\ldots,15,\nonumber
\\
\X_j & = & \eb_j,\qquad j > 15,\nonumber
\end{eqnarray}
where $(\Z_j)_{1}^{3}$ and $(\eb_j)_{1}^{40}$ were i.i.d.\
$N(\zero,\I)$ and $\tau= 0.1$. In this simulation, only coordinates
1 to 5, 6 to 10 and 11 to 15 have nonzero coefficients. These
$x$-variables are uncorrelated across a group, but share the same
correlation within a group. Because the within group correlation is
high, but less than~1, the simulation is ideal for exploring the
effects of strong repressibility.

\begin{figure}[b]

\includegraphics{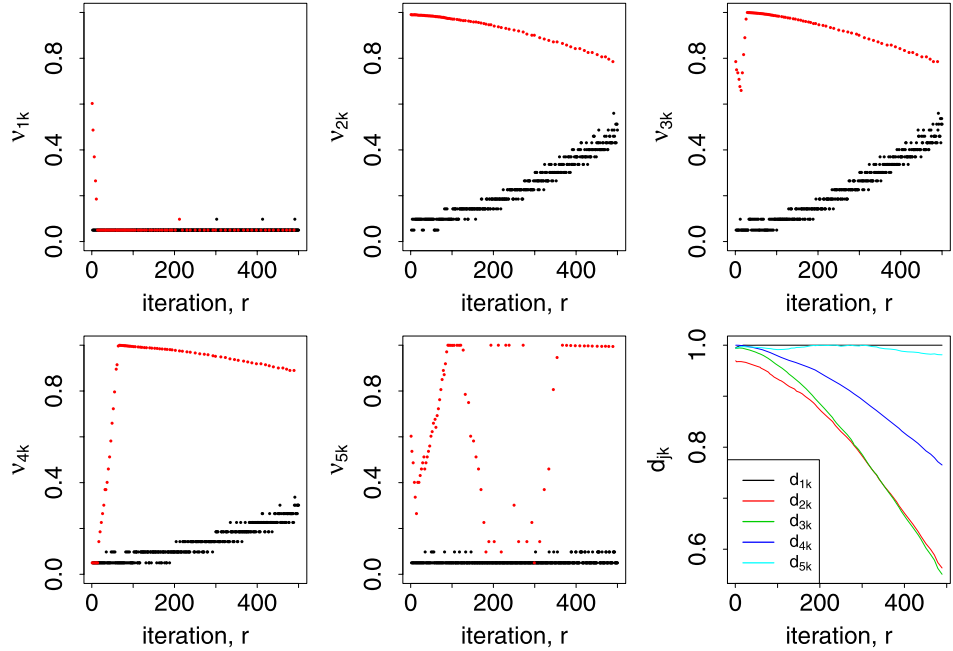}

\caption{First 5 panels display $(\nu_{j,k})_{j=1}^5$ for the
first 5 coefficients from simulation \protect\mref{hastie.sim}: red points
are iterations $r$ where the descent direction
$k\in\{1,\ldots,5\}$. Variables 2 and 3 are never selected due to
their excessively large $\nu_{j,k}$ step sizes: an artifact of the
correlation between the 5 variables. The last panel (bottom right)
displays $(d_{j,k})_{j=1}^5$ for those iterations~$r$ where
$k\in\{1,4,5\}$.}\label{figure5}
\end{figure}

Figure~\ref{figure5} displays results from fitting Algorithm~\ref
{A:L2BoostPath}
for $M^*=500$ iterations with $\nu=0.05$. The first 5
panels are the values $(\nu_{j,k})_{j=1}^5$ against the iteration
$r=1,\ldots,500$, with points colored in red indicating iterations $r$
where $k\in\{1,\ldots,5\}$ and $k$ is used generically to denote the
current descent direction. Notationally, the descent at iteration $r$
is along $k$ for a step size of $\nu_{l,k}$, at which point $l$
becomes more favorable than $k$ and the descent switches to $l$, the
next critical direction. The value plotted, $\nu_{j,k}\le\nu_{l,k}$,
is the step size for $j=1,\ldots,5$.

Whenever the selected coordinate is from the first group of variables
(we are referring to the red points) one of the coordinates $j=1, 4,
5$ achieves a~small $\nu_{j,k}$ value. However, coordinates $j=2$ and
$j=3$ maintain very large values throughout all iterations. This is
despite the fact that the two coordinates generally have large
values of $d_{j,k}$, especially during the early iterations (see the
bottom right panel). This suggests that 1, 4 and 5 become active
variables at some point in the solution path, whereas coordinates 2
and 3 are never selected (indeed, this is exactly what happened). We
can conclude that coordinates 2 and 3 are being strongly repressed by
$k\in\{1,4,5\}$. Interestingly, coordinate 4 also appears to be
repressed at later iterations of the algorithm. Observe how its
$d_{j,k}$ values decrease with increasing $r$ (blue line in bottom
right panel), and that its $\nu_{j,k}$ values are only small at earlier
iterations. Thus, we can also conclude that coordinates $\{1,5\}$
eventually repress coordinate 4 as well.

We note that the number of iterations $M^*=500$ used in the example is
not very large, and if \l2boost were run for a longer period of time,
coordinates 2 and~3 will eventually enter the solution path (panels 2
and~3 of Figure~\ref{figure5} show evidence of this already happening
with $\nu_{j,k}$ steadily decreasing as $r$ increases). However, doing
so leads to overfitting and poor test-set performance (we provide
evidence of this shortly). Using different values of $\nu$ also did
not resolve the problem. Thus, similar to the lasso, we find that
\l2boost is unable to select entire groups of correlated
variables. Like the lasso this means it also will perform suboptimally
in highly correlated settings. In the next section we introduce a
simple way of adding $L_2$-regularization as a way to correct this
deficiency.

\section{Elastic net boosting}\label{S:elasticBoost}

The tendency of the lasso to select only a handful of variables from
among a group of correlated variables was noted
in \citet{Zou:Hast:2005}. To address this
deficiency, \citet{Zou:Hast:2005} described an optimization problem
different from the classical lasso framework. Rather than relying
only on $L_1$-penalization, they included an additional
$L_2$-regularization parameter designed to encourage a ridge-type
grouping effect, and termed the resulting estimator ``the elastic
net.'' Specifically, for a fixed $\lambda>0$ (the ridge parameter)
and a fixed $\lambda_0>0$ (the lasso parameter), the elastic net was
defined as
%
\begin{equation}
\label{enet.def}\qquad
\bhat_\enet=(1+\lambda) \mathop{\arg\min}_{\b\in\RR^p}
\Biggl\{ \Vert\y- \X\b\Vert^2 + \lambda\sum_{k=1}^p\beta_k^2
+ \lambda_0 \sum_{k=1}^p |\beta_k| \Biggr\}.
\end{equation}
To calculate the elastic net, \citet{Zou:Hast:2005} showed
that \mref{enet.def} could be recast as a lasso optimization problem
by replacing the original data with suitably constructed augmented
values. They replaced $\y$ $(n\times1)$ and $\X$ $(n\times p)$ with
augmented values $\ys$ and $\Xs$, defined as follows:
%
\begin{equation}
\label{data.augment}
\ys=
\left[\matrix{
\y
\cr
0
\vspace*{-3pt}\cr
\vdots
\cr
0}
\right]_{(n+p)\times1},\qquad
\Xs= \frac{1}{\sqrt{1+\lambda}}
\left[\matrix{
\X
\cr
\sqrt{\lambda} \I}
\right]_{(n+p)\times p}
= [\Xs_1,\ldots,\Xs_p].\hspace*{-35pt}
\end{equation}
The elastic net optimization can be written in terms of the augmented
data by reparameterizing $\b$ as $\b^*=\b\sqrt{1+\lambda}$.
By Lemma~1 of \citet{Zou:Hast:2005}, it follows that \mref{enet.def}
can be expressed as
\[
\bhat_\enet=\sqrt{1+\lambda} \mathop{\arg\min}_{\b\in\RR
^p}
\Biggl\{ \Vert\y^* - \X^*\b\Vert^2
+ \frac{\lambda_0}{\sqrt{1+\lambda}} \sum_{k=1}^p |\beta_k| \Biggr\},
\]
which is an $L_1$-optimization problem that can be solved using the
lasso.

One explanation for why the elastic net is so successful in correlated
problems is due to its decorrelation property. Let
$R_{j,k}^*=\X_{j}^{*T}\Xs_k$. Because the data is standardized such
that $\X_j^T\X_j=\X_k^T\X_k=1$ [recall \mref{data.standardization}],
we have
\[
R_{j,k}^*
=
\cases{
\displaystyle\frac{\X_j^T\X_k}{1+\lambda}
= \frac{R_{j,k}}{1+\lambda},
&\quad if $j\ne k$,
\vspace*{2pt}\cr
\displaystyle\frac{\X_j^T\X_j + \lambda}{1+\lambda}
= 1,
&\quad if $j = k$.}
\]
One can see that $\lambda$ is a decorrelation parameter,
with larger values reducing the correlation between coordinates.
\citet{Zou:Hast:2005} argued that this effect promotes a ``grouping
property'' for the elastic net that overcomes the lasso's inability to
select groups of correlated variables.

We believe that decorrelation is an important component of the elastic
net's success. However, we will argue that in addition to its role in
decorrelation, $\lambda$ has a surprising connection to repressibility
that further explains its role in regularizing the elastic
net.

The argument for the elastic net follows as a special case (the limit)
of a~generalized \l2boost procedure we refer to as \eboost. The
\eboost\ algorithm is a modification of \l2boost applied to the
augmented problem. To implement \eboost\ one runs \l2boost on the
augmented data \mref{data.augment}, adding a post-processing step to
rescale the coefficient solution path: see Algorithm~\ref{A:eBoost}
for a precise description. For arbitrarily small $\nu$, the solution
path for \eboost\ approximates the elastic net, but for general
$0<\nu\le1$, \eboost\ represents a novel extension of \l2boost. We
study the general \eboost\ algorithm, for arbitrary $0<\nu\le1$, and
present a detailed explanation of how $\lambda$ imposes
$L_2$-regularization.

\begin{algorithm}[t]
\caption{\eboost}\label{A:eBoost}
\begin{algorithmic}[1]
\STATE Augment the data \mref{data.augment}.
Set $F_{0,i}^*=0$ for $i=1,\ldots,n+p$.
\STATE Run Algorithm \ref{A:L2BoostPath} for $M$ iterations using the
augmented data.
\STATE Let $F_{M,i}^*$ denote the $M$-step predictor (discard
$F_{M,i}^*$ for $i>n$).
Let $\beta_{M,k}^*$ denote the $M$-step coefficient estimate.
\STATE Rescale the regression estimates: $\beta_{M,k}=\sqrt{1+\lambda
}\beta_{M,k}^*$.
\end{algorithmic}
\end{algorithm}

\subsection{\texorpdfstring{How $\lambda$ regularizes the solution path}{How lambda regularizes the solution path}}
To study the effect $\lambda$ has on \eboost's solution path we
consider in detail how $\lambda$ effects $m_{j,k}^*$, the number of
steps to favorability [defined as in \mref{critical.step.size} but
with $\y$ and $\X$ replaced by their augmented values $\ys$ and
$\Xs$]. At initialization, the gradient-correlation for $j\ne k$ is
\begin{eqnarray*}
\rho_{j,1}^*
&=& \X_j^{*T}(\ys-\Fhat^{*}_0)
\\
&=& \frac{1}{\sqrt{1+\lambda}}\X_j^T\y- \frac{1}{\sqrt{1+\lambda}}
\Biggl(\sum_{i=1}^n x_{i,j}F_{0,i}^* + \sqrt{\lambda} F_{0,n+j}^*\Biggr).
\end{eqnarray*}
In the special case when $F_{0,i}^*=0$, corresponding to the
first descent of the algorithm,
\[
\rho_{j,1}^*
= \frac{1}{\sqrt{1+\lambda}}\X_j^T\y
= \frac{1}{\sqrt{1+\lambda}}\rho_{j,1}.
\]
Therefore, $d_{j,k}^*=\rho_{j,1}/\rho_{k,1}=d_{j,k}$, and hence
\[
m_{j,k}^*
=\floor\biggl[1+\frac{\log|d_{j,k}-R_{j,k}^*|
-\log(1-R_{j,k}^*\sign(d_{j,k}-R_{j,k}^*) )}{\log(1-\nu)} \biggr].
\]
This equals the number of steps in the original (nonaugmented)
problem but where $\X$ is replaced with variables decorrelated by a
factor of $\sqrt{1+\lambda}$. For large values of $\lambda$ this
addresses the problem seen in Figure~\ref{figure5}.
Recall we argued that
$m_{j,k}$ can became inflated due to the near equality of $d_{j,k}$
with~$R_{j,k}$. However, $R_{j,k}^*=R_{j,k}/\sqrt{1+\lambda}$ shrinks
to zero with increasing $\lambda$, which keeps $m_{j,k}^*$ from
becoming inflated.

This provides one explanation for $\lambda$'s role in regularization,
at least for the case when $\lambda$ is large. But we now suggest
another theory that applies for both small and large $\lambda$. We
argue that regularization is imposed not just by decorrelation, but
through a combination of decorrelation and reversal of repressibility.
Thus $\lambda$'s role is more subtle than our previous argument
suggests.

To show this, let us suppose that near-repressibility holds. We
assume therefore that $R_{j,k}=d_{j,k}(1+\d)$ for some small $|\d|<1$.
Then,
%
\begin{eqnarray}
\label{repress.decorr}
&&\log|d_{j,k}-R_{j,k}^*|
-\log\bigl(1-R_{j,k}^*\sign(d_{j,k}-R_{j,k}^*) \bigr)\nonumber
\\
&& \qquad=
\underbrace{\biggl[\log|d_{j,k}| +
\log\biggl| 1-\frac{1+\d}{\sqrt{1+\lambda}} \biggr| \biggr]}_{\mathrm{Repressibility\
effect}}
\\
&&\hphantom{\qquad=}
\underbrace{- \log\biggl(1-\frac{R_{j,k}}{\sqrt{1+\lambda}}\sign\biggl(R_{j,k}
\biggl[\frac{1}{1+\d}-\frac{1}{\sqrt{1+\lambda}}\biggr] \biggr)\biggr)}_{\mathrm
{Decorrelation\ effect}}.\nonumber
\end{eqnarray}
The first term on the right captures the effect of
repressibility. When $\d$ is small, $\lambda$ plays a crucial
role in controlling its size. If $\lambda=0$, the expression reduces
to $\log|d_{j,k}| + \log|\d|$ which converges to $-\infty$ as
$|\d|\rightarrow0$; thus precluding~$j$ from being selected [keep in
mind that \mref{repress.decorr} is divided by $\log(1-\nu)$, which is
negative; thus $m_{j,k}^*\rightarrow\infty$]. On the other hand, any
$\lambda>0$, even a~relatively small value, ensures that the expression
remains small even for arbitrarily small $\d$, thus reversing the
effect of repressibility.

The second term on the right of \mref{repress.decorr} is related to
decorrelation. If $1+\lambda>(1+\d)^2$ (which holds if $\lambda$ is
large enough when $\d>0$, or for all $\lambda>0$ if $\d<0$), the term
reduces to
\[
- \log\biggl(1-\frac{R_{j,k}}{\sqrt{1+\lambda}}\sign(R_{j,k})\biggr),
\]
which remains bounded when $\lambda>0$ if $R_{j,k}\rightarrow1$.
On the other hand, if $1+\lambda<(1+\d)^2$, the term reduces to
\[
- \log\biggl(1+\frac{R_{j,k}}{\sqrt{1+\lambda}}\sign(R_{j,k})\biggr),
\]
which remains bounded if $R_{j,k}\rightarrow1$ and shrinks in
absolute size as $\lambda$ increases.

Taken together, these arguments show $\lambda$ imposes
$L_2$-regularization through a combination of decorrelation and the
reversal of repressibility which applies even when $\lambda$ is
relatively small.

These arguments apply to the first descent. The general case when
$F_{0,i}^*\ne0$ requires a detailed analysis of $d_{j,k}^*$. In
general,
\[
d_{j,k}^*=
\frac{\X_j^T\y- \sum_{i=1}^n x_{i,j}F_{0,i}^* - \sqrt{\lambda}
F_{0,n+j}^*}
{\X_k^T\y- \sum_{i=1}^n x_{i,k}F_{0,i}^* - \sqrt{\lambda} F_{0,n+k}^*}.
\]
We break up the analysis into two cases depending on the size of
$\lambda$. Suppose first that $\lambda$ is small. Then
\[
d_{j,k}^*\asymp
\frac{\X_j^T\y- \sum_{i=1}^n x_{i,j}F_{0,i}^*}
{\X_k^T\y- \sum_{i=1}^n x_{i,k}F_{0,i}^*},
\]
which is the ratio of gradient correlations based on the original $\X$
without pseudo-data. If $j$ is a promising variable, then $d_{j,k}^*$
will be relatively large, and our\vadjust{\goodbreak} argument from above applies. On the
other hand if $\lambda$ is large, then the third term in the numerator
and the denominator of $d_{j,k}^*$ become the dominating terms and
\[
d_{j,k}^*\asymp\frac{F_{0,n+j}^*}{F_{0,n+k}^*}.
\]
The growth rate of $F_{0,i}^*$ for the pseudo data is $O(\nu)$
for a group of variables that are actively being
explored by the algorithm. Thus $|d_{j,k}^*|\asymp1$ and our
previous argument applies.

\subsection{Illustration}
\begin{figure}[b]
\vspace*{6pt}
\includegraphics{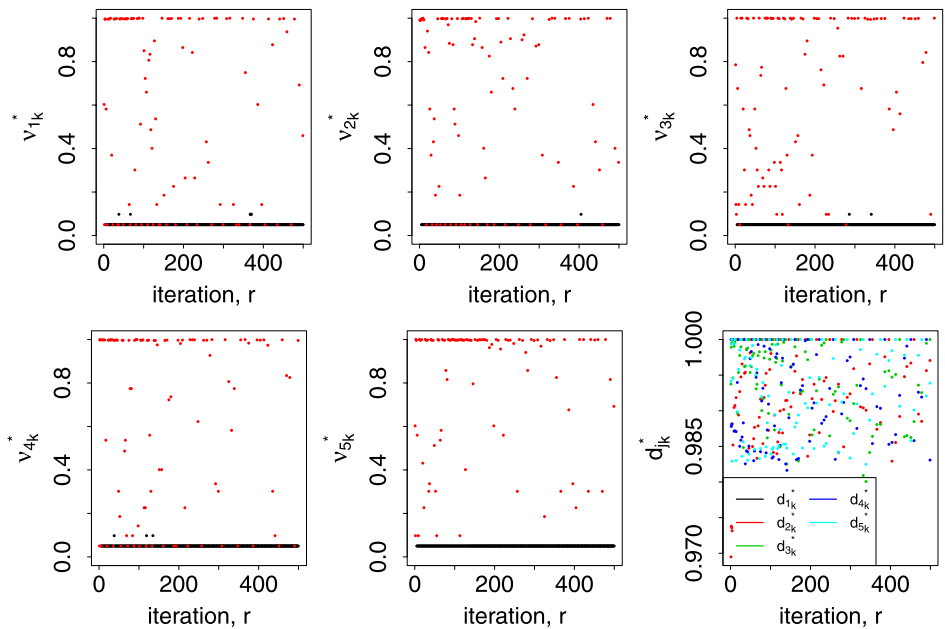}

\caption{\eboost\ applied to simulation \protect\mref{hastie.sim}
(plots are constructed as in Figure~\protect\ref{figure5}). Now each of the
first 5
coordinates are selected and each has $d_{j,k}^*$ values near
one.}\label{figure6}
\end{figure}

As evidence of this, and to demonstrate the effectiveness of \eboost,
we re-analyzed \mref{hastie.sim} using Algorithm~\ref{A:eBoost}. We
used the same parameters as in Figure~\ref{figure5} ($M^*=500$ and
$\nu=0.05$). We set $\lambda=0.5$. The results are displayed in
Figure~\ref{figure6}. In contrast to Figure~\ref{figure5}, notice
that all 5 of the first group of correlated variables achieve small
$\nu_{j,k}^*$ values (and we confirmed that all 5 variables enter the
solution path). It is interesting to note that $d_{j,k}^*$ is nearly
1 for each of these variables.

%

\begin{figure}

\includegraphics{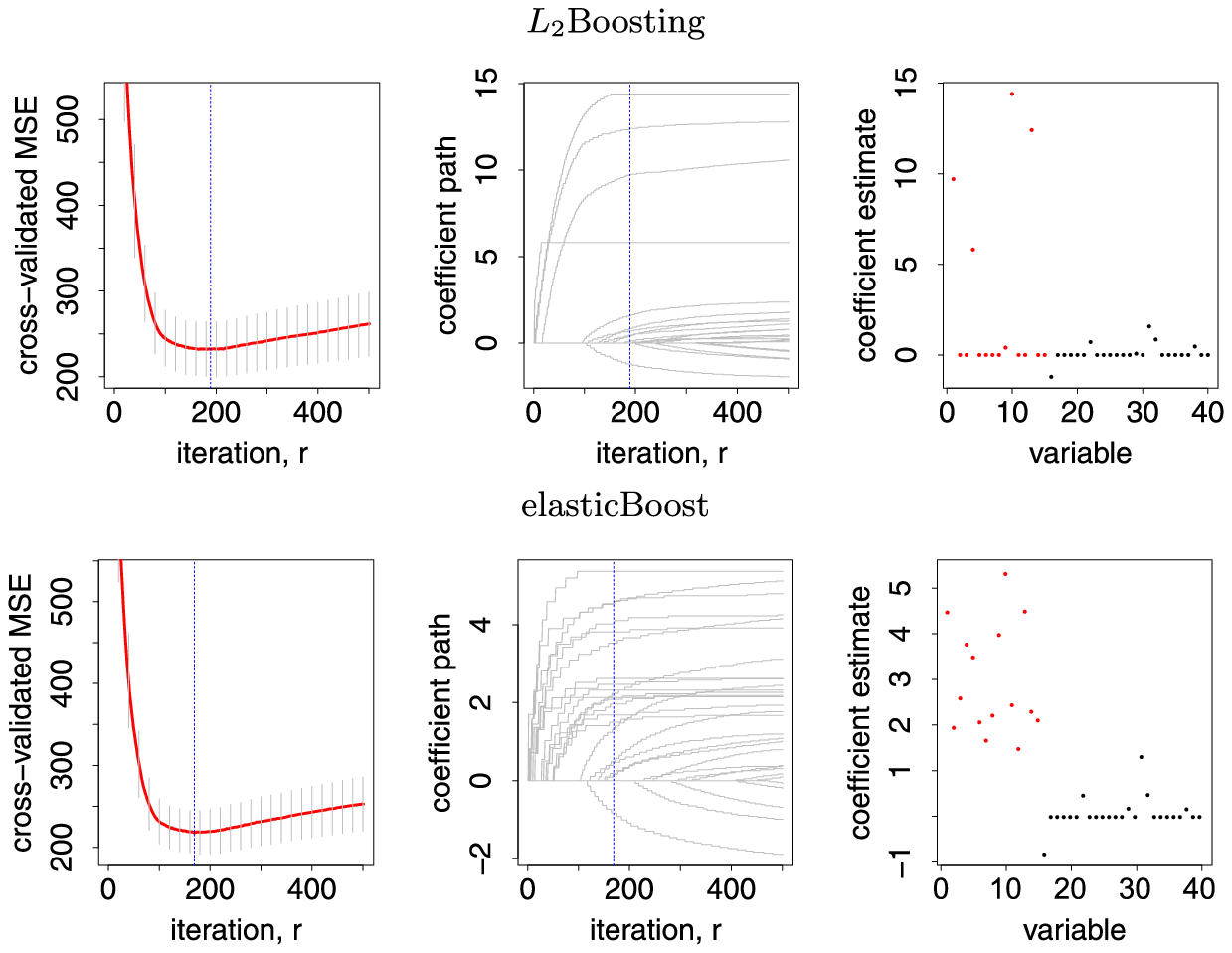}

\caption{\l2boost (top row) versus \eboost\ (bottom row) from
simulation \protect\mref{hastie.sim}.}\label{figure7}
\end{figure}

\begin{figure}

\includegraphics{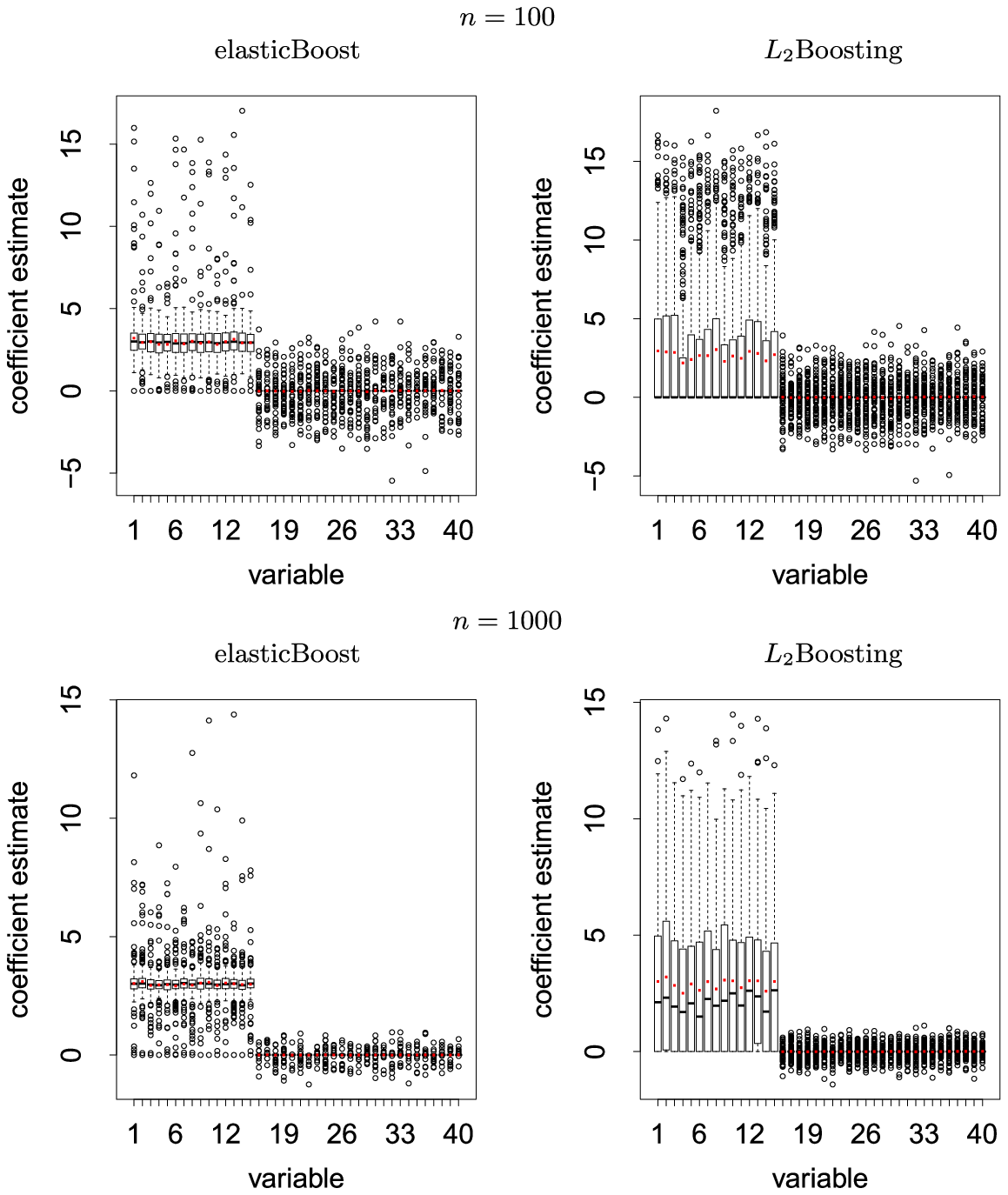}%

\caption{\eboost\ (left) versus \l2boost (right) from
simulation \protect\mref{hastie.sim} for $n=100$ (top) and $n=1000$ (bottom)
based on 250 independent learning samples. The distribution of
coefficient estimates are displayed as boxplots; mean values
are given in red.}\label{figure8}\vspace*{-3pt}
\end{figure}

To compare \l2boost and \eboost\ more evenly, we used 10-fold
cross-validation to determine the optimal number of iterations (for
\eboost, we used doubly-optimized cross-validation to determine both
the optimal number of iterations and the optimal $\lambda$ value;\vadjust{\goodbreak} the
latter was found to equal $\lambda=0.1$). Figure~\ref{figure7}
displays the results. The top row displays \l2boost, while the bottom
row is \eboost\ (fit under the optimized $\lambda$). The minimum
mean-squared-error (MSE) is slightly smaller for \eboost\ (217.9)
than \l2boost (231.7) (first panels in top and bottom rows).
Curiously, the MSE is minimized using about same number of iterations
for both methods (190 for \l2boost and 169 for \eboost). The middle
panels display the coefficient paths. The vertical blue line indicates
the MSE optimized number of iterations. In the case of \l2boost only
4 nonzero coefficients are identified within the optimal number of
steps, whereas \eboost\ finds all 15 nonzero coefficients. This can
be seen more clearly in the right panels which show coefficient
estimates at the optimized stopping time. Not only are all 15
nonzero coefficients identified by \eboost, but their estimated
coefficient values are all roughly near the true value of 3. In
contrast, \l2boost finds only 4 coefficients due to strong
repressibility. Its coefficient estimates are also wildly inaccurate.
While this does not overly degrade prediction error performance (as
evidenced by the first panel), variable selection performance is
seriously impacted.

The entire experiment was then repeated 250 times using 250
independent learning sets. Figure~\ref{figure8} displays the
coefficient estimates from these 250 experiments for \eboost\ (left
side) and \l2boost (right side) as boxplots. The top panel are based
on the original sample size of $n=100$ and the bottom panel use a
larger sample size $n=1000$. The results confirm our previous
finding: \eboost\ is consistently able to group variables and
outperform \l2boost in terms of variable selection.

\begin{figure}

\includegraphics{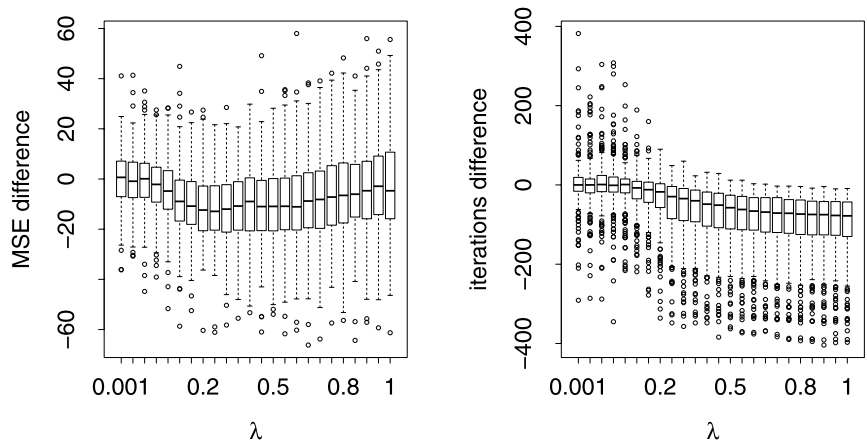}

\caption{Left: difference in test set performance of \l2boost
compared to \eboost. Right: difference in MSE optimized number of
iterations for \l2boost compared to \eboost.}\label{figure9}\vspace*{-3pt}
\end{figure}

Finally, the left panel of Figure~\ref{figure9} displays the
difference in test set MSE for \l2boost and \eboost\ as a
function of $\lambda$ over the 250 experiments ($n=100$).\vadjust{\goodbreak}
Negative values indicate a lower MSE for \eboost, which is generally
the case for larger $\lambda$. The right panel displays the MSE
optimized number of iterations for \l2boost compared to \eboost.
Generally, \eboost\ requires fewer steps as $\lambda$ increases.
This is interesting, because as pointed out, this generally coincides
with better MSE performance.

\section{Discussion}

A key observation is that \l2boost's behavior along a~fixed
descent direction is fully specified with the exception of the descent
length, $L_r$. In Theorem~\ref{criticalpoint.theorem}, we described a
closed form solution for $m_{j,k}$, the number of steps until
favorability, where $k=l_r$ is the currently selected coordinate
direction and $j=l_{r+1}$ is the next most favorable direction.
Theorem~\ref{criticalpoint.theorem} quantifies \l2boost's descent
length, thus allowing us to characterize its solution path as a series
of fixed descents where the next coordinate direction, chosen from all
candidates $j \ne k$, is determined as that with the minimal descent
length $m_{j,k}$ (assuming no ties). Since we choose from among all
directions $j \ne k$, $m_{j,k}$, and equivalently the step length
$\nu_{j,k}$, can be characterized as measures to favorability, a
property of each coordinate at any iteration $r$. These measures are a
function of $\nu$ and the ratio of gradient-correlations $d_{j,k}$ and
the correlation coefficient $R_{j,k}$ relative to the currently
selected direction $k$.

Characterizing the \l2boost solution path by $m_{j,k}$ provides
considerable insight when examining the limiting conditions. When
$m_{j,k} \rightarrow1$, \l2boost exhibits active set cycling, a
property explored in detail in Section~\ref{S:cyclingBehavior}. We
note that this condition is fundamentally a result of the optimization
method which drives $|d_{j,k}|\rightarrow1$ when $\nu$ is arbitrarily
small. This virtually guarantees the notorious slow convergence seen
with infinitesimal forward stagewise algorithms.\vadjust{\goodbreak}

The repressibility condition occurs in the alternative limiting
condition $m_{j,k} \rightarrow\infty$. Repressibility arises when
the gradient correlation ratio $d_{j,k}$ equals the correlation
$R_{j,k}$. When $|R_{j,k}|<1$, $j$ is said to be strongly repressed
by $k$, and while descending along $k$, the absolute
gradient-correlation for $j$ can never be equal to or surpass the
absolute gradient-correlation for~$k$. Strong repressibility plays a
crucial role in correlated settings, hindering variables from being
actively selected. Adding $L_2$ regularization reverses
repressibility and substantially improves variable selection for
\eboost, an \l2boost implementation involving the data augmentation
framework used by the elastic net.\vspace*{-3pt}

\begin{supplement}[id=suppA]
\stitle{Proofs of results from ``Characterizing $L_2$Boosting''\\}
\slink[doi]{10.1214/12-AOS997SUPP} 
\sdatatype{.pdf}
\sfilename{aos997\_supp.pdf}
\sdescription{An online supplementary file contains the detailed
proofs for Theorems~\ref{incremental.operator.theorem} through \ref{dynamic step.two.cycles}. These proofs
make use of various notation described in the paper.\vspace*{-3pt}}
\end{supplement}


\printaddresses

\end{document}